\documentclass[11pt]{article}
\usepackage{amsfonts}
\usepackage{amsfonts}
\usepackage{amsfonts}
\usepackage{mathrsfs}
\usepackage{bm}
\usepackage{amssymb,amsmath} 
\usepackage{cite}
\usepackage{cite}
\usepackage[colorlinks,linkcolor=blue,citecolor=blue]{hyperref}
 \allowdisplaybreaks

\catcode`!=11
\let\!int\int \def\int{\displaystyle\!int}
\let\!lim\lim \def\lim{\displaystyle\!lim}
\let\!sum\sum \def\sum{\displaystyle\!sum}
\let\!sup\sup \def\sup{\displaystyle\!sup}
\let\!inf\inf \def\inf{\displaystyle\!inf}
\let\!cap\cap \def\cap{\displaystyle\!cap}
\let\!max\max \def\max{\displaystyle\!max}
\let\!min\min \def\min{\displaystyle\!min}
\let\!frac\frac \def\frac{\displaystyle\!frac}
\catcode`!=12

\let\oldsection\section
\renewcommand\section{\setcounter{equation}{0}\oldsection}

\allowdisplaybreaks
\def\pf{\it{Proof.}\rm\quad}

\def\R{\mathbb{R}}
\def\N{\mathbb{N}}\def\Z{\mathbb{Z}}

\newtheorem{thm}{Theorem}[section]
\newtheorem{lem}[thm]{Lemma}
\newtheorem{cor}[thm]{Corollary}

\begin{document}
\title {\bf Some Evaluations of Parametric Euler Type Sums of Harmonic Numbers}
\author{
{Junjie Quan\thanks{Email: qjunjiexmu@xujc.com}\quad Ce Xu$^{}$\thanks{Email: cexu2020@ahnu.edu.cn}\quad Xixi Zhang\thanks{Email: XixiZhang2019@163.com}}\\[1mm]
\small * School of Information Science and Technology,  Xiamen University Tan Kah Kee College\\
\small Xiamen
 Fujian 363105, P.R. China\\
 \small $\dag,\ddag$ School of Mathematics and Statistics, Anhui Normal University,\\ \small Wuhu 241002, P.R. China}

\date{}
\maketitle \noindent{\bf Abstract } We establish some identities of Euler related sums. By using these identities, we discuss the closed form representations of sums of harmonic numbers and reciprocal parametric binomial coefficients through parametric harmonic numbers, shifted harmonic numbers and Riemann zeta function with positive integer arguments. In particular we investigate products of quadratic and cubic harmonic numbers and reciprocal parametric binomial coefficients. Some illustrative special cases as well
as immediate consequences of the main results are also considered.
\\[2mm]
\noindent{\bf Keywords} Euler sums; harmonic numbers; binomial coefficients; polylogarithms.
\\[2mm]
\noindent{\bf AMS Subject Classifications (2010):}  05A10, 05A19, 11B65, 11B83, 11M06, 33D60, 33C20.

\section{Introduction}

Let $\R$ and $\mathbb{C}$ denote respectively, the sets of real and complex numbers and let $\N=\{1,2,3,\ldots\}$ be the set of natural numbers, and $\N_0:=\N\cup \{0\}$. The classical harmonic numbers of order $s$ are given by (\cite{A1995})
\begin{align}\label{a1}
H_n^{\left( s \right)}: = \sum\limits_{j = 1}^n {\frac{1}
{{{j^s}}}} \;\;\left( {n \in \N,\;s \in \mathbb{C}} \right)
\end{align}
and
\begin{align}\label{a2}
{H_n}: = H_n^{\left( 1 \right)} = \sum\limits_{j = 1}^n {\frac{1}
{j}}  = \int\limits_0^1 {\frac{{1 - {t^n}}}
{{1 - t}}dt}  = \gamma  + \psi \left( {n + 1} \right)\ \ (n\in \N),
\end{align}
where $\gamma$ denotes the Euler-Mascheroni constant, defined by
\[\gamma  := \mathop {\lim }\limits_{n \to \infty } \left( {\sum\limits_{k = 1}^n {\frac{1}{k}}  - \ln n} \right) =  - \psi \left( 1 \right) \approx {\rm{ 0 }}{\rm{. 577215664901532860606512 }}...,\]
and $\psi \left( z \right)$ denotes the digamma function (or called Psi function ) which is defined as the logarithmic derivative of the well known gamma function:
\begin{align}\label{a3}
\psi \left( z \right) := \frac{d}{{dz}}\left( {\ln \Gamma \left( z \right)} \right) = \frac{{\Gamma '\left( z \right)}}{{\Gamma \left( z \right)}}.
\end{align}

In general, the polygamma function of order $m$ is defined as the $(m+1)$-th derivative of the logarithm of the gamma function:
\begin{align}\label{a4}
{\psi ^{\left( m \right)}}\left( z \right): = \frac{{{d^m}}}
{{d{z^m}}}\psi \left( z \right) = \frac{{{d^{m + 1}}}}
{{d{z^{m + 1}}}}\ln \left( {\Gamma \left( z \right)} \right)\ \ (m\in \N_0; z\in \mathbb{C}\setminus \Z_0^-:=\{0,-1,-2,\ldots\}).
\end{align}
For any $s\in \mathbb{C}\setminus \{0\}$, we assume that
\[{H_0}: = 0,\;H_0^{\left( s \right)}: = 0,\;H_0^{\left( 0 \right)}: = 1.\]
Equation (1.1) can be written in the following form:
\[H_n^{\left( s \right)} = \zeta \left( s \right) - \zeta \left( {s,n + 1} \right)\;\;\left( {{\Re} \left( s \right) > 1;n \in \N} \right)\]
by recalling the well-known (easily derivable) relationship between the Riemann Zeta function $\zeta(s)$ (or called Euler-Riemann zeta function) and the Hurwitz (or
generalized) Zeta function $\zeta(s,a)$ (see \cite{A2000})
\begin{align}\label{a5}
\zeta \left( s \right) = \zeta \left( {s,n + 1} \right) + \sum\limits_{j = 1}^n {\frac{1}
{{{j^s}}}}\ \ (n\in \N_0).
\end{align}
The Riemann Zeta function or Euler-Riemann zeta function (for more details, see for instance, \cite{A1972,A2000,SC2012})
\[\zeta (s): = \sum\limits_{n = 1}^\infty  {\frac{1}
{{{n^s}}}} \;\;\left( {{\Re} \left( s \right) > 1} \right)\]
is probably the most important, fascinating, challenging and mysterious object of modern mathematics, in spite of its utter simplicity. This function is defined over the complex plane and plays a pivotal role in analytic number theory having applications in physics, probability theory, applied statistics and other fields of mathematics. There is an enormous amount of literature on the Riemann zeta function.

A well-known (and potentially useful) relationship between the polygamma functions ${\psi ^{\left( m \right)}}\left( z \right)$ and the generalized Zeta function (or Hurwitz Zeta function) $\zeta(s,a)$ is given by
\begin{align}\label{a6}
{\psi ^{\left( m \right)}}\left( z \right) = {\left( { - 1} \right)^{m + 1}}m!\sum\limits_{n = 0}^\infty  {\frac{1}
{{{{\left( {n + z} \right)}^{m + 1}}}}}  = {\left( { - 1} \right)^{m + 1}}m!\zeta \left( {m + 1,z} \right).
\end{align}
It is also easy to have the following expression
\[{\psi ^{\left( m \right)}}\left( {n + 1} \right) - {\psi ^{\left( m \right)}}\left( 1 \right) = {\left( { - 1} \right)^{m }}m!H_n^{\left( {m + 1} \right)}\ \ (m,n\in \N_0),\]
which immediately yields another expression for $H^{(m)}_n$ as follows
\begin{align}\label{a7}
H_n^{\left( m \right)} = \frac{{{{\left( { - 1} \right)}^{m-1}}}}
{{\left( {m - 1} \right)!}}\left\{ {{\psi ^{\left( {m - 1} \right)}}\left( {n + 1} \right) - {\psi ^{\left( {m - 1} \right)}}\left( 1 \right)} \right\}\ \ (m\in \N;n\in \N_0).
\end{align}
Replacing $n$ by $\alpha\ (\alpha\in \mathbb{C}\setminus \Z^-;\Z^-:=\{-1,-2,\ldots\})$ in above equations (1.2) and (1.7), we can define the generalized $\alpha$-th shifted harmonic number of order $m$ by
\begin{align}\label{a8}
H_\alpha ^{\left( m \right)}: = \frac{{{{\left( { - 1} \right)}^{m - 1}}}}{{\left( {m - 1} \right)!}}\left( {{\psi ^{\left( {m - 1} \right)}}\left( {\alpha  + 1} \right) - {\psi ^{\left( {m - 1} \right)}}\left( 1 \right)} \right)\quad (2 \le m \in \N)
\end{align}
and
\begin{align}\label{a9}
{H_\alpha }\equiv H_\alpha ^{\left( 1 \right)}: = \psi \left( {\alpha  + 1} \right) + \gamma=\sum\limits_{k = 1}^\infty  {\left( {\frac{1}{k} - \frac{1}{{k + \alpha }}} \right)}.
\end{align}
It is easy to otain the following expression by using the relation (1.6) and the definition of Hurwitz Zeta function $\zeta(s,a+1)$
\begin{align}\label{a10}
H_\alpha ^{\left( m \right)} = \zeta \left( m \right) - \zeta \left( {m,\alpha  + 1} \right),\;\;H_{n + \alpha }^{\left( m \right)} = H_\alpha ^{\left( m \right)} + \sum\limits_{j = 1}^n {\frac{1}
{{{{\left( {j + \alpha } \right)}^m}}}},\;2 \le m \in \N .
\end{align}
The generalized harmonic functions (also called parametric harmonic numbers or partial sums of Hurwitz Zeta function) are defined by (see \cite{C2011,RS2002,S2011,W2016})
\begin{align}\label{a11}
H_n^{\left( s\right)}\left( a \right): = \sum\limits_{j = 1}^n {\frac{1}
{{{{\left( {j + a} \right)}^s}}}}\ (n\in \N;s\in \mathbb{C};a\in \mathbb{C}\setminus\Z^-),
\end{align}
so that, obviously, $H_n^{\left( s \right)}\left( 0 \right) = H_n^{\left( s \right)}.$
By using the equation (1.6),  we may write the generalized harmonic numbers, $H_n^{\left( m\right)}\left( a \right)$ , in terms of polygamma functions
\begin{align}\label{a12}
H_n^{\left( m \right)}\left( a \right) = \frac{{{{\left( { - 1} \right)}^{m - 1}}}}
{{\left( {m - 1} \right)!}}\left\{ {{\psi ^{\left( {m - 1} \right)}}\left( {n + a + 1} \right) - {\psi ^{\left( {m - 1} \right)}}\left( {a + 1} \right)} \right\}\ (m\in\N;n\in \N_0).
\end{align}

A generalized binomial coefficient $\left( {\begin{array}{*{20}{c}}
   a  \\
   b  \\
\end{array}} \right)\;\left( {a,b \in \mathbb{C}} \right)$ is defined, in terms of the familiar (Euler's) gamma function, by
\begin{align}\label{a13}
\left( {\begin{array}{*{20}{c}}
   a  \\
   b  \\
\end{array}} \right): = \frac{{\Gamma \left( {a + 1} \right)}}{{\Gamma \left( {b + 1} \right)\Gamma \left( {a - b + 1} \right)}},\;\;a,b \in \mathbb{C},
\end{align}
which, in the special case when $b=n,n\in \N_0$, yields
\[\left( {\begin{array}{*{20}{c}}
   a  \\
   0  \\
\end{array}} \right): = 1,\left( {\begin{array}{*{20}{c}}
   a  \\
   n  \\
\end{array}} \right): = \frac{{a\left( {a - 1} \right) \cdots \left( {a - n + 1} \right)}}{{n!}},\;n \in \N.\]

Here we are interested in evaluating closed form and integral representations of Euler type sums containing both harmonic numbers and binomial coefficients. In this paper we will develop identities, closed form representations of harmonic numbers and reciprocal binomial coefficients of the form:
\begin{align}\label{a14}
W_k^{\left( {a,b} \right)}\left( {{m_1},{m_2}, \cdots ,{m_r},p} \right): = \sum\limits_{n = 1}^\infty  {\frac{{H_n^{\left( {{m_1}} \right)}H_n^{\left( {{m_2}} \right)} \cdots H_n^{\left( {{m_r}} \right)}}}
{{{{\left( {n + a} \right)}^p}\left( {\begin{array}{*{20}{c}}
   {n + k + b}  \\
   k  \\
 \end{array} } \right)}}}\ \ (a,b\in \mathbb{C}\setminus\Z^-)
 \end{align}
for $r=1,p\in \{0,1\}$ and $r=2,m_i=1,p\in \{0,1\}$, where $p\in \N_0, {k,r,{m_i} \in \N,i = 1,2, \ldots ,r}$ and $p+k>1$.
When $a=b=0$, then the parametric Euler type sum $W_k^{\left( {a,b} \right)}\left( {{m_1},{m_2}, \cdots ,{m_r},p} \right)$ reduces to the classical Euler type sum $W_k\left( {{m_1},{m_2}, \cdots ,{m_r},p} \right)$ which is defined by (see \cite{XMZ2016})
\begin{align}\label{a15}
{W_k}\left( {{m_1},{m_2}, \cdots ,{m_r},p} \right) := \sum\limits_{n = 1}^\infty  {\frac{{H_n^{\left( {{m_1}} \right)}H_n^{\left( {{m_2}} \right)} \cdots H_n^{\left( {{m_r}} \right)}}}{{{n^p}\left( {\begin{array}{*{20}{c}}
   {n + k}  \\
   k  \\
\end{array}} \right)}}} ,\;\left(p\in \N_0, {k,r,{m_i} \in \N} \right)
 \end{align}
with $p+k>1$.
There has recently been renewed interest in the study of series involving binomial coefficients and a number of authors have obtained either closed form representation or integral representation for some particular cases of these series. The interested reader is referred to \cite{C2011,CS2011,C2013,KCDC2021, JCC2008,BG1996, S2011,S2015,SD2012,SS2011,SC2012,WX2019,XMZ2016}, and references therein. For results on alternating quadratic harmonic number sums see \cite{S2015}. When $k=0$ in (1.15), then the Euler type sum ${W_k}\left( {{m_1},{m_2}, \cdots ,{m_r},p} \right)$ reduces to the classical Euler sum $S_{{\bf m},p}\ \ ({\bf m}:=\{m_1,m_2,\ldots,m_r\})$ defined by (see \cite{FS1998})
\begin{align}\label{a16}
S_{{\bf m},p}:=\sum\limits_{n = 1}^\infty  {\frac{{H_n^{\left( {{m_1}} \right)}H_n^{\left( {{m_2}} \right)} \cdots H_n^{\left( {{m_r}} \right)}}}
{{{n^p}}}}.
 \end{align}
There are many works investigating Euler sums $S_{{\bf m},p}$ and related variants involving harmonic numbers, see for example \cite{DB2015,FS1998,F2005,M2014,W2017,X2016}, and references therein. There are fewer results for sums of the type (1.14). In this paper, we will prove that all Euler type sums
\[W_k^{\left( {a,b} \right)}\left( {1,p} \right),\;W_k^{\left( {a,a} \right)}\left( {m,p} \right),\;W_k^{\left( {a,a} \right)}\left( {1,1,p} \right)\]
for $p=0,1$  can be expressed as a rational linear combination of products of zeta values, shifted harmonic numbers and  parametric harmonic numbers.

\section{Some Lemmas}

The following lemmas will be useful in the development of the main theorems.

\begin{lem}(\cite[Thm. 4.1]{CX2017}) For $m\in \N,\ a,b\in \mathbb{C}\setminus\Z^-$ and $x\in (-1,1)$. Then the following identity holds:
\begin{align}\label{b1}
\int\limits_0^x {{H_m}\left( {t,a} \right)} {t^{n + b - 1}}dt =& \sum\limits_{k = 1}^{m - 1} {\frac{{{{\left( { - 1} \right)}^{k - 1}}{x^{n + b}}}}{{{{\left( {n + b} \right)}^k}}}} {H_{m + 1 - k}}\left( {x,a} \right)\nonumber\\
& + \frac{{{{\left( { - 1} \right)}^{m - 1}}}}{{{{\left( {n + b} \right)}^m}}}\left\{ {{x^{n + b}}{H_1}\left( {x,a} \right) + \sum\limits_{k = 1}^n {\frac{{{x^{k + a + b}}}}{{k + a + b}}}  - {H_1}\left( {x,a + b} \right)} \right\},
\end{align}
where the function ${H_m}\left( x,a \right)$ is defined by
\begin{align}\label{b2}
{H_s}\left( x,a \right) := \sum\limits_{n = 1}^\infty  {\frac{{{x^{n + a}}}}{{{{\left( {n + a} \right)}^s}}}} \quad (\Re(s) > 1,\;a\notin \Z^-
,\;x \in \left[ { - 1,1} \right])
\end{align}
with ${H_1}\left( {x,a} \right): = \sum\limits_{n = 1}^\infty  {\frac{{{x^{n + a}}}}{{n + a}}}$ and $x\in [-1,1)$.
\end{lem}

\begin{cor}\label{corb-2} For $m\in \N,\ b\in \mathbb{C}\setminus\Z^-$ and $x\in (-1,1)$, we have
\begin{align}\label{bb1}
\int\limits_0^x {{}{{\rm Li}_m}\left( t \right){t^{n + b - 1}}dt}  =& \sum\limits_{i = 1}^{m - 1} {\frac{{{{\left( { - 1} \right)}^{i - 1}}}}
{{{{\left( {n + b} \right)}^i}}}{x^{n + b}}{{\rm Li}_{m + 1 - i}}\left( x \right)}  + \frac{{{{\left( { - 1} \right)}^{m - 1}}}}
{{{{\left( {n + b} \right)}^m}}}\sum\limits_{j = 1}^n {\frac{{{x^{j + b}}}}
{{j + b}}} \nonumber\\
& + \frac{{{{\left( { - 1} \right)}^{m - 1}}}}
{{{{\left( {n + b} \right)}^m}}}\left\{ {{x^{n + b}}{{\rm Li}_1}\left( x \right) - {H_1}\left( {x,b} \right)} \right\},
\end{align}
where ${\rm Li}_s(x)$ denotes the polylogarithm function which is defined by
\[{\rm Li}{_s}\left( x \right) = \sum\limits_{n = 1}^\infty  {\frac{{{x^n}}}{{{n^s}}}} ,\;\Re (s)\geq1, - 1 \le x < 1,\]
with ${\rm Li}_1(x)=-\ln(1-x)$. Obviously, ${{\rm Li}_s}\left( x \right) = {H_s}\left( {x,0} \right).$
\end{cor}
\pf Setting $a=0$ in \eqref{b1}, the result is \eqref{bb1}.\hfill$\square$

\begin{lem}\label{lemb-3} (\cite{CX2016,CX2017}) Let $s$ be positive integer and $a\notin \Z^-$, then we have
\begin{align}\label{b3}
\sum\limits_{n = 1}^\infty  {\frac{{{y^n}\sum\limits_{j = 1}^{n - 1} {\frac{{{x^{n - j}}}}{j}}  + {x^n}\sum\limits_{j = 1}^{n - 1} {\frac{{{y^{n - j}}}}{j}} }}{{{{\left( {n + a} \right)}^s}}}}  =& s{\rm Li}{_{s + 1}}\left( {a,xy} \right) - \sum\limits_{j = 1}^s {{\rm Li}{_j}\left( {a,x} \right){\rm Li}{_{s + 1 - j}}\left( {a,y} \right)}\nonumber
\\&  + {\rm Li}{_s}\left( {a,xy} \right)\left( {{\rm Li}{_1}\left( x \right) + {\rm Li}{_1}\left( y \right)} \right),
\end{align}
where $x,y \in \left[ { - 1,1} \right)$ and the parametric polylogarithm function ${{\rm Li}{_s}\left( {a,x} \right)}$ is defined by
\[{\rm Li}{_s}\left( {a,x} \right) = \sum\limits_{n = 1}^\infty  {\frac{{{x^n}}}{{{{\left( {n + a} \right)}^s}}}} ,\;\Re (s)\geq1, - 1 \le x < 1.\]
\end{lem}
\pf Taking $x=y=q$ and $b=a$ in Theorem 2.2 in the reference \cite{CX2016}, then letting $q$ approach 1, we obtain formula \eqref{b3}. \hfill$\square$

Putting $x=y$ in above equation, we arrive at the conclusion that ($s\geq 2$)
\begin{align}\label{b4}
\sum\limits_{n = 1}^\infty  {\frac{{{x^n}}}{{{{\left( {n + a} \right)}^s}}}} \sum\limits_{j = 1}^{n - 1} {\frac{{{x^{n - j}}}}{j}}
 &=\frac{s}{2}{\rm Li}{_{s + 1}}\left( {a,{x^2}} \right) + {\rm Li}{_s}\left( {a,{x^2}} \right){\rm Li}{_1}\left( x \right) - {\rm Li}{_s}\left( {a,x} \right){\rm Li}{_1}\left( {a,x} \right)
\nonumber \\
           &\quad - \frac{1}{2}\sum\limits_{j = 2}^{s - 1} {{\rm Li}{_j}\left( {a,x} \right){\rm Li}{_{s + 1 - j}}\left( {a,x} \right)}.
\end{align}
Note that when $s>1$, $\mathop {\lim }\limits_{x \to 1}{\rm{L}}{{\rm{i}}_s}\left( {a,1} \right) = \zeta \left( {s,a + 1} \right)$. When $x$ approach 1, we arrive at the conclusion that
\begin{align}\label{b5}
\mathop {\lim }\limits_{x \to 1} \left\{ {{\rm{L}}{{\rm{i}}_s}\left( {a,{x^2}} \right){\rm{L}}{{\rm{i}}_1}\left( x \right) - {\rm{L}}{{\rm{i}}_s}\left( {a,{x}} \right){\rm{L}}{{\rm{i}}_1}\left( {a,x} \right)} \right\} = \zeta \left( {s,a + 1} \right)H_a.
\end{align}
Hence, letting $x\rightarrow 1$ in \eqref{b4}, we obtain
\begin{align}\label{b6}
\sum\limits_{n = 1}^\infty  {\frac{{{H_n}}}{{{{\left( {n + a} \right)}^s}}}}
 &=\frac{s}{2}\zeta \left( {s + 1,a + 1} \right) - \frac{1}{2}\sum\limits_{j = 1}^{s - 2} {\zeta \left( {s - j,a + 1} \right)} \zeta \left( {j + 1,a + 1} \right)
\nonumber \\
           &\quad + \zeta \left( {s,a + 1} \right)H_a  + \sum\limits_{n = 1}^\infty  {\frac{1}{{n{{\left( {n + a} \right)}^s}}}}.
\end{align}
By using partial fraction decomposition, we get
\begin{align}\label{b7}
\sum\limits_{n = 1}^\infty  {\frac{1}
{{n{{\left( {n + a} \right)}^s}}}}  = \frac{{{H_a}}}
{{{a^s}}} - \sum\limits_{j = 2}^s {\frac{{\zeta \left( {j,a + 1} \right)}}
{{{a^{s + 1 - j}}}}} \quad (s\in \N;a \in \mathbb{C}\setminus\Z_0^-).
\end{align}
Therefore, we know that the parametric Euler sums \[\sum\limits_{n = 1}^\infty  {\frac{{{H_n}}}{{{{\left( {n + a} \right)}^s}}}}\quad (2\leq s\in \N)\] can be expressed in term of parametric harmonic numbers and shifted harmonic numbers.

\begin{lem}\label{lemb-4} For any real $x\in (-1,1)$, then the following identity holds:
\begin{align}\label{b8}
\sum\limits_{n = 1}^\infty  {{H_n}H_n^{\left( 2 \right)}{x^n}}  = \frac{1}{{1 - x}}\left\{ {2{\rm Li}{_3}\left( x \right) - \ln \left( {1 - x} \right){\rm Li}{_2}\left( x \right) - \sum\limits_{n = 1}^\infty  {\frac{{{H_n}}}{{{n^2}}}{x^n}} } \right\}.
\end{align}
\end{lem}
\pf In \cite{SX2017}, Xin etc. proved the result
\begin{align}\label{b9}
\sum\limits_{n = 1}^\infty  {{H_n}{H^{(m)} _n}{x^n}}  = \frac{1}{{1 - x}}\left\{ {\sum\limits_{n = 1}^\infty  {\frac{{{H_n}}}{{{n^m}}}{x^n}}  - \sum\limits_{n = 1}^\infty  {\frac{1}{{{n^m}}}\left( {\sum\limits_{k = 1}^n {\frac{{{x^k}}}{k}} } \right)}  - \zeta \left( m \right)\ln \left( {1 - x} \right)} \right\}.
\end{align}
Then, we consider the nested sum
\[\sum\limits_{n = 1}^\infty  {\frac{{{y^n}}}{{{n^m}}}\left( {\sum\limits_{k = 1}^n {\frac{{{x^k}}}{{{k^p}}}} } \right)} ,\;x,y \in \left[ { - 1,1} \right),\;m,p \in \N.\]
By taking the sum over complementary pairs of summation indices, we obtain a simple reflection formula
\begin{align}\label{b10}
\sum\limits_{n = 1}^\infty  {\frac{{{y^n}}}{{{n^m}}}\left( {\sum\limits_{k = 1}^n {\frac{{{x^k}}}{{{k^p}}}} } \right)}  + \sum\limits_{n = 1}^\infty  {\frac{{{x^n}}}{{{n^p}}}\left( {\sum\limits_{k = 1}^n {\frac{{{y^k}}}{{{k^m}}}} } \right)}  = {\rm Li}{_p}\left( x \right){\rm Li}{_m}\left( y \right) + {\rm Li}{_{p + m}}\left( {xy} \right).
\end{align}
Setting $p=1,m=2,y=1$ in above equation we get
\begin{align}\label{b11}
\sum\limits_{n = 1}^\infty  {\frac{1}{{{n^2}}}\left( {\sum\limits_{k = 1}^n {\frac{{{x^k}}}{k}} } \right)}  + \sum\limits_{n = 1}^\infty  {\frac{{H_n^{\left( 2 \right)}}}{n}{x^n}}  =  - \ln \left( {1 - x} \right)\zeta \left( 2 \right) + {\rm{L}}{{\rm{i}}_3}\left( x \right).
\end{align}
On the other hand, by the definition of polylogarithm function and Cauchy product of power series, we have
\begin{align}\label{b12}
\sum\limits_{n = 1}^\infty  {\frac{{H_n^{\left( 2 \right)}}}{n}{x^n}} & =\int\limits_0^x {\frac{{{\rm{L}}{{\rm{i}}_2}\left( t \right)}}{{t\left( {1 - t} \right)}}dt}  = \int\limits_0^x {\frac{{{\rm{L}}{{\rm{i}}_2}\left( t \right)}}{t}dt}  + \int\limits_0^x {\frac{{{\rm{L}}{{\rm{i}}_2}\left( t \right)}}{{1 - t}}dt} \nonumber\\
& = {\rm{L}}{{\rm{i}}_3}\left( x \right) - \ln \left( {1 - x} \right){\rm{L}}{{\rm{i}}_2}\left( x \right) - \int\limits_0^x {\frac{{{{\ln }^2}\left( {1 - t} \right)}}{t}} dt\nonumber\\
& = 3{\rm{L}}{{\rm{i}}_3}\left( x \right) - \ln \left( {1 - x} \right){\rm{L}}{{\rm{i}}_2}\left( x \right) - 2\sum\limits_{n = 1}^\infty  {\frac{{{H_n}}}{{{n^2}}}{x^n}} .
\end{align}
Then, substituting \eqref{b12} into \eqref{b11} yields
\begin{align}\label{b13}
\sum\limits_{n = 1}^\infty  {\frac{1}{{{n^2}}}\left( {\sum\limits_{k = 1}^n {\frac{{{x^k}}}{k}} } \right)}  = 2\sum\limits_{n = 1}^\infty  {\frac{{{H_n}}}{{{n^2}}}{x^n}}  + \ln \left( {1 - x} \right){\rm{L}}{{\rm{i}}_2}\left( x \right) - 2{\rm{L}}{{\rm{i}}_3}\left( x \right) - \ln \left( {1 - x} \right)\zeta \left( 2 \right).
\end{align}
Hence, taking $m=2$ in \eqref{b9} and combining formula \eqref{b13} we may deduce the desired result.
The proof of Lemma \ref{lemb-4} is thus completed.\hfill$\square$

\begin{lem}(\cite[Thm. 2.1]{X2017}) For $m\in \N$ and $a>0$, we have the recurrence relation
\begin{align}\label{b14}
{Y_m}\left( a \right) = \left( {m - 1} \right)!\sum\limits_{i = 0}^{m - 1} {\frac{{{Y_i}\left( a \right)}}
{{i!}}H_a^{\left( {m - i} \right)}} ,\;{Y_0}\left( a \right) = 1,
\end{align}
where ${Y_m}\left( a \right)$ is defined by the following integral
\begin{align}\label{b15}
\frac{{{Y_m}\left( a \right)}}{a}:={\left( { - 1} \right)^m}\int\limits_0^1 {{x^{a - 1}}{{\ln }^m}\left( {1 - x} \right)dx}.
\end{align}
\end{lem}

From the recurrence relation \eqref{b14}, we obtain the following identities
\begin{align}
&\int\limits_0^1 {{x^{a - 1}}\ln \left( {1 - x} \right)dx}  =  - \frac{{{H_a}}}{a}, \label{b16}\\
&\int\limits_0^1 {{x^{a - 1}}{{\ln }^2}\left( {1 - x} \right)dx}  = \frac{{H_a^2 + H_a^{\left( 2 \right)}}}
{a},\label{b17}\\
 &\int\limits_0^1 {{x^{a - 1}}{{\ln }^3}\left( {1 - x} \right)dx}  =  - \frac{{H_a^3 + 3{H_a}H_a^{\left( 2 \right)} + 2H_a^{\left( 3 \right)}}}
{a}, \label{b18}\\
&\int\limits_0^1 {{x^{a - 1}}{{\ln }^4}\left( {1 - x} \right)dx}  = \frac{{H_a^4 + 6H_a^2H_a^{\left( 2 \right)} + 8{H_a}H_a^{\left( 3 \right)} + 3{{\left( {H_a^{\left( 2 \right)}} \right)}^2} + 6H_a^{\left( 4 \right)}}}
{a}.\label{b19}
 \end{align}

\section{Identities for Euler related sums}
We now prove the following theorems.

\begin{thm} For $a,b \in \mathbb{C}\setminus\Z_0^-$ with $a\neq b$. Then the following identity holds:
\begin{align}\label{c1}\sum\limits_{n = 1}^\infty  {\frac{{{H_n}}}
{{\left( {n + a} \right)\left( {n + b} \right)}}}  = \frac{1}
{{b - a}}\left( {\frac{{{H_a}}}
{a} - \frac{{{H_b}}}
{b}} \right) + \frac{{H_a^2 - H_b^2}}
{{2\left( {b - a} \right)}} + \frac{{\zeta \left( {2,a + 1} \right) - \zeta \left( {2,b + 1} \right)}}
{{2\left( {b - a} \right)}}.
\end{align}
\end{thm}
\pf First we note that we may rewrite the series on the left hand side of \eqref{c1} as
\begin{align}\label{c2}\sum\limits_{n = 1}^\infty  {\frac{{{H_n}}}
{{\left( {n + a} \right)\left( {n + b} \right)}}}  = \frac{1}
{{b - a}}\mathop {\lim }\limits_{x \to 1} \left\{ {\sum\limits_{n = 1}^\infty  {\left( {\frac{{{x^n}}}
{{n + a}} - \frac{{{x^n}}}
{{n + b}}} \right)\left( {\sum\limits_{j = 1}^n {\frac{{{x^{n - j}}}}
{j}} } \right)} } \right\}.
\end{align}
On the other hand, from Lemma \ref{lemb-3}, we deduce that
\begin{align}\label{c3}
\sum\limits_{n = 1}^\infty  {\frac{{{x^n}}}
{{n + a}}} \sum\limits_{j = 1}^{n - 1} {\frac{{{x^{n - j}}}}
{j}}={\rm Li}_1(a,x^2){\rm Li}_1(x)+\frac{1}{2}{\rm Li}_2(a,x^2)-\frac{1}{2}({\rm Li}_1(a,x))^2.
\end{align}
Substituting \eqref{c3} into \eqref{c2}, by a simple calculation, the result is
\begin{align}\label{c4}
\sum\limits_{n = 1}^\infty  {\frac{{{H_n}}}
{{\left( {n + a} \right)\left( {n + b} \right)}}}  =& \sum\limits_{n = 1}^\infty  {\frac{1}
{{n\left( {n + a} \right)\left( {n + b} \right)}}}  + \frac{1}
{{2\left( {b - a} \right)}}\sum\limits_{n = 1}^\infty  {\left\{ {\frac{1}
{{{{\left( {n + a} \right)}^2}}} - \frac{1}
{{{{\left( {n + b} \right)}^2}}}} \right\}}\nonumber \\
& + \frac{1}
{2}\left( {\sum\limits_{n = 1}^\infty  {\frac{1}
{{\left( {n + a} \right)\left( {n + b} \right)}}} } \right)\left( {a\sum\limits_{n = 1}^\infty  {\frac{1}
{{n\left( {n + a} \right)}} + b\sum\limits_{n = 1}^\infty  {\frac{1}
{{n\left( {n + b} \right)}}} } } \right).
\end{align}
Using the definitions of shifted harmonic number and Hurwitz Zeta function we obtain the desired result. \hfill$\square$
\begin{thm} For $m,k\in \N$ and $a\in \mathbb{C}\setminus\Z^-$. Then
\begin{align}\label{c5}
\sum\limits_{n = 1}^\infty  {\frac{{H_n^{\left( m \right)}}}
{{\left( {n + a} \right)\left( {n + a + k} \right)}}}  = \frac{1}
{k}\left\{ \begin{gathered}
  \sum\limits_{n = 1}^\infty  {\frac{1}
{{{n^m}\left( {n + a} \right)}}}  + \sum\limits_{j = 1}^{m - 1} {{{\left( { - 1} \right)}^{j - 1}}\zeta \left( {m + 1 - j} \right)H_{k - 1}^{\left( j \right)}\left( a \right)}  \hfill \\
   + {\left( { - 1} \right)^{m - 1}}{H_a}H_{k - 1}^{\left( m \right)}\left( a \right) + {\left( { - 1} \right)^{m - 1}}\sum\limits_{i = 1}^{k - 1} {\frac{{H_i^{\left( 1 \right)}\left( a \right)}}
{{{{\left( {i + a} \right)}^m}}}}  \hfill \\
\end{gathered}  \right\},
\end{align}
when $a\in \mathbb{C}\setminus\Z_0^-$, we have (see \cite[Eq. (1.11)]{X2017})
\begin{align}\label{c6}
\int\limits_0^1 {{x^{a - 1}}{{\rm Li}_m}\left( x \right)dx} =\sum\limits_{n = 1}^\infty  {\frac{1}
{{{n^m}\left( {n + a} \right)}}}  = \sum\limits_{l = 1}^{m - 1} {{{\left( { - 1} \right)}^{l - 1}}\frac{{\zeta \left( {m + 1 - l} \right)}}
{{{a^l}}}}  + {\left( { - 1} \right)^{m - 1}}\frac{{{H_a}}}
{{{a^m}}}.
\end{align}
\end{thm}
\pf By the definition of polylogarithm function and Cauchy product formula, we can verify that
\begin{align}\label{c7}
\frac{{{\rm Li}{_m}\left( t \right)}}{{1 - t}} = \sum\limits_{n = 1}^\infty  {{H^{(m)} _n}{t^n}}, \ t\in (-1,1).
\end{align}
Multiplying \eqref{c7} by $t^{a-1}-t^{a+k-1}$ and integrating over $(0,1)$, we obtain
\begin{align}\label{c8}
k\sum\limits_{n = 1}^\infty  {\frac{{H_n^{\left( m \right)}}}
{{\left( {n + a} \right)\left( {n + a + k} \right)}}}  = \int\limits_0^1 {{x^{a - 1}}{\text{L}}{{\text{i}}_m}\left( x \right)dx}  + \sum\limits_{i = 1}^{k - 1} {\int\limits_0^1 {{x^{a + i - 1}}{\text{L}}{{\text{i}}_m}\left( x \right)dx} } .
\end{align}
In Corollary \ref{corb-2}, letting $x\rightarrow 1$ and noting the fact that
\[\mathop {\lim }\limits_{x \to 1} \left\{ {{x^{n + a}}{\text{L}}{{\text{i}}_1}\left( x \right) - {H_1}\left( {x,a} \right)} \right\} = {H_a},\]
we arrive at the conclusion that
\begin{align}\label{c9}
\int\limits_0^1 {{x^{a + n - 1}}{\text{L}}{{\text{i}}_m}\left( x \right)dx}  = \sum\limits_{j = 1}^{m - 1} {\frac{{{{\left( { - 1} \right)}^{j - 1}}}}
{{{{\left( {n + a} \right)}^j}}}\zeta \left( {m + 1 - j} \right)}  + {\left( { - 1} \right)^{m - 1}}\frac{{{H_a}}}
{{{{\left( {n + a} \right)}^m}}} + {\left( { - 1} \right)^{m - 1}}\frac{{H_n^{\left( 1 \right)}\left( a \right)}}
{{{{\left( {n + a} \right)}^m}}}.
\end{align}
Combining \eqref{c6}, \eqref{c8} and \eqref{c9}, the result is \eqref{c5}.\hfill$\square$

From Theorem 2.3, we can get the following corollary.
\begin{cor}\label{corc-3} For integers $m,k,r\in \N$ and $k\neq r$, we have
\begin{align}\label{c10}
\sum\limits_{n = 1}^\infty  {\frac{{{H^{(m)} _n}}}{{n\left( {n + k} \right)}}}  = \frac{1}{k}\left\{ {\zeta \left( {m + 1} \right) + \sum\limits_{j = 1}^{m - 1} {{{\left( { - 1} \right)}^{j - 1}}\zeta \left( {m + 1 - j} \right){H^{(j)} _{k - 1}}}  + {{\left( { - 1} \right)}^{m - 1}}\sum\limits_{i = 1}^{k - 1} {\frac{{{H_i}}}{{{i^m}}}} } \right\}.
\end{align}
\begin{align}\label{c11}
\sum\limits_{n = 1}^\infty  {\frac{{{H^{(m)} _n}}}{{\left( {n + r} \right)\left( {n + k} \right)}}}  = \frac{1}{{k - r}}\left\{ \begin{array}{l}
 \sum\limits_{j = 1}^{m - 1} {{{\left( { - 1} \right)}^{j - 1}}\zeta \left( {m + 1 - j} \right)\left( {{H^{(j)} _{k - 1}} - {H^{(j)} _{r - 1}}} \right)}  \\
  + {\left( { - 1} \right)^{m - 1}}\left( {\sum\limits_{i = 1}^{k - 1} {\frac{{{H_i}}}{{{i^m}}}}  - \sum\limits_{i = 1}^{r - 1} {\frac{{{H_i}}}{{{i^m}}}} } \right) \\
\end{array} \right\}.
\end{align}
\end{cor}

Putting $m=1,2$ in \eqref{c5}, we give the following identities
\begin{align}\label{c12}
\sum\limits_{n = 1}^\infty  {\frac{{{H_n}}}
{{\left( {n + a} \right)\left( {n + a + k} \right)}}}  = \frac{1}
{k}\left\{ {{H_a}H_k^{\left( m \right)}\left( {a - 1} \right) - \frac{{H_k^{\left( 1 \right)}\left( a \right)}}
{{k + a}} + \frac{{{{\left( {H_k^{\left( 1 \right)}\left( a \right)} \right)}^2} + H_k^{\left( 2 \right)}\left( a \right)}}
{2}} \right\},
\end{align}
\begin{align}\label{c13}
\sum\limits_{n = 1}^\infty  {\frac{{H_n^{\left( 2 \right)}}}
{{\left( {n + a} \right)\left( {n + a + k} \right)}}}  = \frac{1}
{k}\left\{ {\zeta \left( 2 \right)H_k^{\left( 1 \right)}\left( {a - 1} \right) - {H_a}H_k^{\left( 2 \right)}\left( {a - 1} \right) - \sum\limits_{i = 1}^{k - 1} {\frac{{H_i^{\left( 1 \right)}\left( a \right)}}
{{{{\left( {i + a} \right)}^2}}}} } \right\}.
\end{align}

\begin{thm}\label{thmc-4} For $k\in \N$ and $a>0$, we have
\begin{align}\label{c14}
\sum\limits_{n = 1}^\infty  {\frac{{H_n^2}}
{{\left( {n + a} \right)\left( {n + a + k} \right)}}}  = \frac{1}
{k}\left\{ \begin{gathered}
  \zeta \left( 2 \right)H_k^{\left( 1 \right)}\left( {a - 1} \right) - {H_a}H_k^{\left( 2 \right)}\left( {a - 1} \right) \hfill \\
   - \sum\limits_{i = 1}^{k - 1} {\frac{{H_i^{\left( 1 \right)}\left( a \right)}}
{{{{\left( {i + a} \right)}^2}}}}  + \sum\limits_{j = 1}^k {\frac{{H_{a + j - 1}^2 + H_{a + j - 1}^{\left( 2 \right)}}}
{{a + j - 1}}}  \hfill \\
\end{gathered}  \right\}.
\end{align}
\end{thm}
\pf From \cite{L1974}, we know that the generating function of Stirling numbers of the first kind is given by
\begin{align}\label{c15}
{\ln ^{m + 1}}\left( {1 - x} \right) = {\left( { - 1} \right)^{m + 1}}\left( {m + 1} \right)!\sum\limits_{n = m + 1}^\infty  {s\left( {n,m + 1} \right)} \frac{{{x^n}}}{{n!}},m\in \N_0,x\in [-1,1),
\end{align}
where ${s\left( {n,k} \right)}$ denotes the (unsigned) Stirling number of the first kind (see \cite{L1974}), and we have
\begin{align*}
& s\left( {n,1} \right) = \left( {n - 1} \right)!,s\left( {n,2} \right) = \left( {n - 1} \right)!{H_{n - 1}},s\left( {n,3} \right) = \frac{{\left( {n - 1} \right)!}}{2}\left[ {H_{n - 1}^2 - {H^{(2)} _{n - 1}}} \right],\\
&s\left( {n,4} \right) = \frac{{\left( {n - 1} \right)!}}{6}\left[ {H_{n - 1}^3 - 3{H_{n - 1}}{H^{(2)} _{n - 1}} + 2{H^{(3)} _{n - 1}}} \right], \\
&s\left( {n,5} \right) = \frac{{\left( {n - 1} \right)!}}{{24}}\left[ {H_{n - 1}^4 - 6{H^{(4)} _{n - 1}} - 6H_{n - 1}^2{H^{(2)} _{n - 1}} + 3(H^{(2)} _{n - 1})^2 + 8H_{n - 1}^{}{H^{(3)} _{n - 1}}} \right].
\end{align*}
Hence, we obtain the following formula
\begin{align}\label{c16}
{\ln ^3}\left( {1 - x} \right) = 3\sum\limits_{n = 1}^\infty  {\frac{{H_n^{\left( 2 \right)} - H_n^2}}
{{n + 1}}{x^{n + 1}}}.
\end{align}
Differentiating this equality above, we deduce that
\begin{align}\label{c17}
\frac{{{{\ln }^2}\left( {1 - x} \right)}}
{{1 - x}} = \sum\limits_{n = 1}^\infty  {\left( {H_n^2 - H_n^{\left( 2 \right)}} \right){x^n}} .
\end{align}
Multiplying \eqref{c7} by $x^{a-1}-x^{a+k-1}$ and integrating over $(0,1)$, yields
\begin{align}\label{c18}
k\sum\limits_{n = 1}^\infty  {\frac{{H_n^2 - H_n^{\left( 2 \right)}}}
{{\left( {n + a} \right)\left( {n + a + k} \right)}}}  = \sum\limits_{j = 1}^k {\int\limits_0^1 {{x^{a + j - 2}}{{\ln }^2}\left( {1 - x} \right)dx} } .
\end{align}
Combining \eqref{b17} with \eqref{c18}, we get
\begin{align}\label{c19}
\sum\limits_{n = 1}^\infty  {\frac{{H_n^2 - H_n^{\left( 2 \right)}}}
{{\left( {n + a} \right)\left( {n + a + k} \right)}}}  = \frac{1}
{k}\sum\limits_{j = 1}^k {\frac{{H_{a + j - 1}^2 + H_{a + j - 1}^{\left( 2 \right)}}}
{{a + j - 1}}},
\end{align}
which together with \eqref{c13} gives \eqref{c14} and finishes the proof of Theorem \ref{thmc-4}.\hfill$\square$
\begin{thm}\label{thmc-5} For $k\in \N$ and $a>0$, then the following identities hold:
\begin{align}\label{c20}
\sum\limits_{n = 1}^\infty  {\frac{{{H_n}H_n^{\left( 2 \right)}}}{{\left( {n + a} \right)\left( {n + a + k} \right)}}}  = \frac{1}{k}\left\{ \begin{array}{l}
 \sum\limits_{i = 0}^{k - 1} {\frac{1}{{a + i}}\sum\limits_{n = 1}^\infty  {\frac{{{H_{n + a + i}}}}{{{n^2}}}} }  - \sum\limits_{i = 0}^{k - 1} {\frac{{H_{a + i}^2 + H_{a + i}^{\left( 2 \right)}}}{{2{{\left( {a + i} \right)}^2}}}}  \\
  - \zeta \left( 2 \right)H_k^{\left( 2 \right)}\left( {a - 1} \right) + \sum\limits_{i = 0}^{k - 1} {\frac{{{H_{a + i}}}}{{{{\left( {a + i} \right)}^3}}}}  \\
 \end{array} \right\},
\end{align}
\begin{align}\label{c21}
\sum\limits_{n = 1}^\infty  {\frac{{H_n^3}}{{\left( {n + a} \right)\left( {n + a + k} \right)}}}  = \frac{1}{k}\left\{ \begin{array}{l}
 \sum\limits_{i = 1}^k {\frac{{H_{a + i - 1}^3 + 3{H_{a + i - 1}}H_{a + i - 1}^{\left( 2 \right)} + 2H_{a + i - 1}^{\left( 3 \right)}}}{{a + i - 1}}}  \\
  + 3\sum\limits_{i = 0}^{k - 1} {\frac{1}{{a + i}}\sum\limits_{n = 1}^\infty  {\frac{{{H_{n + a + i}}}}{{{n^2}}}} }  + 3\sum\limits_{i = 1}^{k - 1} {\frac{{{H_{a + i}}}}{{{{\left( {a + i} \right)}^3}}}}  \\
  - \frac{3}{2}\sum\limits_{i = 0}^{k - 1} {\frac{{H_{a + i}^2 + H_{a + i}^{\left( 2 \right)}}}{{{{\left( {a + i} \right)}^2}}}}  - \zeta \left( 2 \right)H_{k - 1}^{\left( 2 \right)}\left( a \right) \\
  - 2\sum\limits_{i = 1}^{k - 1} {\frac{{H_i^{\left( 1 \right)}\left( a \right)}}{{{{\left( {i + a} \right)}^3}}}}  - 2\zeta \left( 3 \right)H_{k - 1}^{\left( 1 \right)}\left( a \right) \\
  - 2{H_a}H_{k - 1}^{\left( 3 \right)}\left( a \right) - 2\frac{{\zeta \left( 3 \right)}}{a} - \frac{{\zeta \left( 2 \right)}}{{{a^2}}} + \frac{{{H_a}}}{{{a^3}}} \\
 \end{array} \right\}.
\end{align}
\end{thm}
\pf First, we are ready to prove the formula \eqref{c20}. Multiplying \eqref{b8} by $x^{a-1}-x^{a+k-1}$ and integrating over the interval $(0,1)$, we can arrive at the conclusion that
\begin{align}\label{c22}
 k\sum\limits_{n = 1}^\infty  {\frac{{{H_n}H_n^{\left( 2 \right)}}}{{\left( {n + a} \right)\left( {n + a + k} \right)}}} =& 2\sum\limits_{i = 0}^{k - 1} {\int\limits_0^1 {{x^{a + i - 1}}{\rm{L}}{{\rm{i}}_3}\left( x \right)dx} }  - \sum\limits_{i = 0}^{k - 1} {\sum\limits_{n = 1}^\infty  {\frac{{{H_n}}}{{{n^2}\left( {n + a + i} \right)}}} } \nonumber \\
  &- \sum\limits_{i = 0}^{k - 1} {\int\limits_0^1 {{x^{a + i - 1}}{\rm{ln}}\left( {1 - x} \right){\rm{L}}{{\rm{i}}_2}\left( x \right)dx} } .
  \end{align}
From \eqref{b16}, \eqref{b17} and \eqref{c6}, we readily find that
\begin{align}\label{c23}
\int\limits_0^1 {{x^{a - 1}}{\rm{L}}{{\rm{i}}_3}\left( x \right)dx}  = \frac{{\zeta \left( 3 \right)}}{a} - \frac{{\zeta \left( 2 \right)}}{{{a^2}}} + \frac{{{H_a}}}{{{a^3}}},
 \end{align}
\begin{align}\label{c24}
 \sum\limits_{n = 1}^\infty  {\frac{{{H_{n + a}}}}{{n\left( {n + a} \right)}}}  =&  - \sum\limits_{n = 1}^\infty  {\frac{1}{n}} \int\limits_0^1 {{x^{n + a - 1}}{\rm{ln}}\left( {1 - x} \right)dx} = \int\limits_0^1 {{x^{a - 1}}{\rm{l}}{{\rm{n}}^2}\left( {1 - x} \right)dx}\nonumber \\
  =& \frac{{H_a^2 + H_a^{\left( 2 \right)}}}{a}.
\end{align}
Therefore, by using equations \eqref{b16} and \eqref{c24}, we deduce that
\begin{align}\label{c25}
\int\limits_0^1 {{x^{\alpha  - 1}}\ln \left( {1 - x} \right){\rm{L}}{{\rm{i}}_2}\left( x \right)dx}  = & - \sum\limits_{n = 1}^\infty  {\frac{{{H_{n + a }}}}{{{n^2}\left( {n + a } \right)}}}
 = \frac{1}{a}\sum\limits_{n = 1}^\infty  {\frac{{{H_{n + \alpha }}}}{{n\left( {n + a } \right)}}}  - \frac{1}{a }\sum\limits_{n = 1}^\infty  {\frac{{{H_{n + a }}}}{{{n^2}}}}\nonumber \\
 =&\frac{{H_a^2 + H_a^{\left( 2 \right)}}}{a}- \frac{1}{a }\sum\limits_{n = 1}^\infty  {\frac{{{H_{n + a }}}}{{{n^2}}}}.
\end{align}
Moreover, using \eqref{c4}, by a direct calculation, then
\begin{align}\label{c26}
\sum\limits_{n = 1}^\infty  {\frac{{{H_n}}}{{{n^2}\left( {n + a} \right)}}}  = \frac{{2\zeta \left( 3 \right)}}{a} - \frac{{\zeta \left( 2 \right)}}{{{a^2}}} + \frac{{{H_a}}}{{{a^3}}} - \frac{{H_a^2 + H_a^{\left( 2 \right)}}}{{2{a^2}}}.
\end{align}
Hence, substituting \eqref{c23}, \eqref{c25} and \eqref{c26} into \eqref{c22}, after simplification, we may easily deduce the result \eqref{c20}.

Next, we prove the formula \eqref{c21}. By a similar argument as in the proof of \eqref{c19}, it is easily shown that
\begin{align}\label{c27}\left( {p - 1} \right)!\sum\limits_{n = p - 1}^\infty  {\frac{{s\left( {n + 1,p} \right)}}{{n!\left( {n + a} \right)\left( {n + a + k} \right)}}}  = \frac{1}{k}\sum\limits_{i = 1}^k {\frac{{{Y_{p - 1}}\left( {a + i - 1} \right)}}{{a + i - 1}}} \;\left( {a > 0} \right).
\end{align}
Setting $p=4$ in above equation we obtain
\begin{align}\label{c28}\sum\limits_{n = 1}^\infty  {\frac{{H_n^3 - 3{H_n}H_n^{\left( 2 \right)} + 2H_n^{\left( 3 \right)}}}{{\left( {n + a} \right)\left( {n + a + k} \right)}} = \frac{1}{k}} \sum\limits_{i = 1}^k {\frac{{H_{a + i - 1}^3 + 3{H_{a + i - 1}}H_{a + i - 1}^{\left( 2 \right)} + 2H_{a + i - 1}^{\left( 3 \right)}}}{{a + i - 1}}} .
\end{align}
Letting $m=3$ in \eqref{c5}, we can get
\begin{align}\label{c29}\sum\limits_{n = 1}^\infty  {\frac{{H_n^{\left( 3 \right)}}}{{\left( {n + a} \right)\left( {n + a + k} \right)}} = \frac{1}{k}} \left\{ \begin{array}{l}
 \frac{{\zeta \left( 3 \right)}}{a} - \frac{{\zeta \left( 2 \right)}}{{{a^2}}} + \frac{{{H_a}}}{{{a^3}}} + \zeta \left( 3 \right)H_{k - 1}^{\left( 1 \right)}\left( a \right) \\
  - \zeta \left( 2 \right)H_{k - 1}^{\left( 2 \right)}\left( a \right) + {H_a}H_{k - 1}^{\left( 3 \right)}\left( a \right) + \sum\limits_{i = 1}^{k - 1} {\frac{{H_i^{\left( 1 \right)}\left( a \right)}}{{{{\left( {i + a} \right)}^3}}}}  \\
 \end{array} \right\}.
 \end{align}
Thus, combining \eqref{c20}, \eqref{c28} and \eqref{c29} we have the result \eqref{c21}. The proof of Theorem \ref{thmc-5} is finished.\hfill$\square$\\
\section{The main results}
In this section, we give some closed form sums of ${W^{(a,b)}_k}\left( {{m_1},{m_2}, \cdots ,{m_r},p} \right)$
through parametric harmonic numbers, shifted harmonic numbers and zeta values. First, we consider the expansion
\begin{align}\label{d1}\frac{1}{{\prod\limits_{i = 1}^k {\left( {n + {a_i}} \right)} }} = \sum\limits_{j = 1}^k {\frac{{{A_j}}}{{n + {a_j}}}}\ \ (k\in \N_0; a_i \in \mathbb{C}\setminus\Z^-) \end{align}
where
\begin{align}\label{d2}{A_j} = \mathop {\lim }\limits_{n \to  - {a_j}} \frac{{n + {a_j}}}{{\prod\limits_{i = 1}^k {\left( {n + {a_i}} \right)} }} = \prod\limits_{i = 1,i \ne j}^k {{{\left( {{a_i} - {a_j}} \right)}^{ - 1}}}.\end{align}
Taking $a_i=a+i$ in \eqref{d2} we obtain
\begin{align}\label{d3}
{A_r} = {\left( { - 1} \right)^{r + 1}}\frac{r}
{{k!}}\left( {\begin{array}{*{20}{c}}
   k  \\
   r \\
 \end{array} } \right),\end{align}
\begin{align}\label{d4}\frac{1}
{{\prod\limits_{i = 1}^k {\left( {n + a + i} \right)} }} = \sum\limits_{r = 1}^k {{{\left( { - 1} \right)}^{r + 1}}\frac{r}
{{k!}}\left( {\begin{array}{*{20}{c}}
   k  \\
   r  \\
 \end{array} } \right)\frac{1}
{{n + a + r}}} .\end{align}
Furthermore, by using the equation \eqref{d4} and the definition of binomial coefficient, we have the following expansions
\begin{align}
&\frac{1}
{{\left( {\begin{array}{*{20}{c}}
   {n + k + a}  \\
   k  \\
\end{array} } \right)}} = \sum\limits_{r = 1}^k {{{\left( { - 1} \right)}^{r + 1}}r\left( {\begin{array}{*{20}{c}}
   k  \\
   r  \\
\end{array} } \right)\frac{1}
{{n + a + r}}} \ \ (k\in \N_0; a \in \mathbb{C}\setminus\Z^-),\label{d5}\\
&\frac{1}
{{\left( {\begin{array}{*{20}{c}}
   {n + k + a}  \\
   k  \\
 \end{array} } \right)}} = k\sum\limits_{r = 1}^{k - 1} {{{\left( { - 1} \right)}^{r + 1}}r\left( {\begin{array}{*{20}{c}}
   {k - 1}  \\
   r  \\
\end{array} } \right)\frac{1}
{{\left( {n + a + 1} \right)\left( {n + r + 1 + a} \right)}}}\ \ (k\in \N; a \in \mathbb{C}\setminus\Z^-).\label{d6}
\end{align}
Therefore, we can obtain the following relations
\begin{align}\label{d7}
\sum\limits_{n = 1}^\infty  {\frac{{f\left( n \right)}}
{{\left( {\begin{array}{*{20}{c}}
   {n + k + a}  \\
   k  \\
 \end{array} } \right)}}}  = k\sum\limits_{r = 1}^{k - 1} {{{\left( { - 1} \right)}^{r + 1}}r\left( {\begin{array}{*{20}{c}}
   {k - 1}  \\
   r  \\
 \end{array} } \right)\frac{{f\left( n \right)}}
{{\left( {n + a + 1} \right)\left( {n + r + 1 + a} \right)}}} ,
\end{align}
\begin{align}\label{d8}
\sum\limits_{n = 1}^\infty  {\frac{{f\left( n \right)}}
{{{{\left( {n + a} \right)}^p}\left( {\begin{array}{*{20}{c}}
   {n + k + b}  \\
   k  \\
 \end{array} } \right)}}}  = \sum\limits_{r = 1}^k {{{\left( { - 1} \right)}^{r + 1}}r\left( {\begin{array}{*{20}{c}}
   k  \\
   r  \\
\end{array} } \right)\sum\limits_{n = 1}^\infty  {\frac{{f\left( n \right)}}
{{{{\left( {n + a} \right)}^p}\left( {n + b + r} \right)}}} } .
\end{align}
The main result of this paper is embodied in the following theorems.
\begin{thm}\label{thmd-1} For $a,b \in \mathbb{C}\setminus\Z_0^-$ and $p,k\in \N$ with $a\neq b+r\ (r=1,2,\ldots,k)$. Then
\begin{align}\label{d9}
W_k^{\left( {a,b} \right)}\left( {1,p} \right) &= \sum\limits_{n = 1}^\infty  {\frac{{{H_n}}}
{{{{\left( {n + a} \right)}^p}\left( {\begin{array}{*{20}{c}}
   {n + k + b}  \\
   k  \\
 \end{array} } \right)}}} \nonumber\\
 &= \sum\limits_{r = 1}^k {{{\left( { - 1} \right)}^{r + 1}}r\left( {\begin{array}{*{20}{c}}
   k  \\
   r  \\
 \end{array} } \right)\frac{1}
{{{{\left( {a - b - r} \right)}^{p - 1}}}}\left\{ \begin{gathered}
  \frac{1}
{{b + r - a}}\left( {\frac{{{H_a}}}
{a} - \frac{{{H_{b + r}}}}
{{b + r}}} \right) + \frac{{H_a^2 - H_{b + r}^2}}
{{2\left( {b + r - a} \right)}} \hfill \\
   + \frac{{\zeta \left( {2,a + 1} \right) - \zeta \left( {2,b + r + 1} \right)}}
{{2\left( {b + r - a} \right)}} \hfill \\
\end{gathered}  \right\}}\nonumber \\
&\quad - \sum\limits_{r = 1}^k {{{\left( { - 1} \right)}^{r + 1}}r\left( {\begin{array}{*{20}{c}}
   k  \\
   r  \\
 \end{array} } \right)\sum\limits_{j = 2}^p {\frac{1}
{{{{\left( {a - b - r} \right)}^{p + 1 - j}}}}\sum\limits_{n = 1}^\infty  {\frac{{{H_n}}}
{{{{\left( {n + a} \right)}^j}}}} } }.
\end{align}
\end{thm}
\pf Putting $f(n)=H_n$ in \eqref{d8} and using the following partial fraction decomposition
\begin{align}\label{d10}
\frac{1}
{{{{\left( {n + a} \right)}^p}\left( {n + b} \right)}} = \frac{1}
{{{{\left( {a - b} \right)}^{p - 1}}}}\frac{1}
{{\left( {n + a} \right)\left( {n + b} \right)}} - \sum\limits_{j = 2}^p {\frac{1}
{{{{\left( {a - b} \right)}^{p + 1 - j}}}}\frac{1}
{{{{\left( {n + a} \right)}^j}}}},
\end{align}
then combining \eqref{c1}, we may easily deduce the result.\hfill$\square$
\begin{thm}\label{thmd-2} For $a\in \mathbb{C}\setminus\Z^-$ and $k\in \N\setminus\{1\}$. Then the following identity holds:
\begin{align}\label{d11}W_k^{\left( {a,b} \right)}\left( {m,0} \right) = k\sum\limits_{r = 1}^{k - 1} {{{\left( { - 1} \right)}^{r + 1}}\left( {\begin{array}{*{20}{c}}
   {k - 1}  \\
   r  \\
 \end{array} } \right)\left\{ \begin{gathered}
\sum_{i=1}^{m-1} \frac{(-1)^{i-1}}{(a+1)^i}\zeta(m+1-i)-(-1)^m \frac{H_{a+1}}{(a+1)^m}  \hfill \\
   + {\left( { - 1} \right)^{m - 1}}\sum\limits_{i = 1}^{r - 1} {\frac{{H_i^{\left( 1 \right)}\left( {a + 1} \right)}}
{{{{\left( {i + a + 1} \right)}^m}}}}  \hfill \\
   + {\left( { - 1} \right)^{m - 1}}{H_{a + 1}}H_{r - 1}^{\left( m \right)}\left( {a + 1} \right) \hfill \\
   + \sum\limits_{j = 1}^{m - 1} {{{\left( { - 1} \right)}^{j - 1}}\zeta \left( {m + 1 - j} \right)H_{r - 1}^{\left( j \right)}\left( {a + 1} \right)}  \hfill \\
\end{gathered}  \right\}}.
\end{align}
\end{thm}
\pf Taking $f(n)=H^{(m)}_n$ in \eqref{d7}, we have
\begin{align}\label{d12}
W_k^{\left( {a,b} \right)}\left( {m,0} \right) =& \sum\limits_{n = 1}^\infty  {\frac{{H_n^{\left( m \right)}}}
{{\left( {\begin{array}{*{20}{c}}
   {n + k + a}  \\
   k  \\
 \end{array} } \right)}}}\nonumber \\
  =& k\sum\limits_{r = 1}^{k - 1} {{{\left( { - 1} \right)}^{r + 1}}r\left( {\begin{array}{*{20}{c}}
   {k - 1}  \\
   r  \\
 \end{array} } \right)} \sum\limits_{n = 1}^\infty  {\frac{{H_n^{\left( m \right)}}}
{{\left( {n + a + 1} \right)\left( {n + r + 1 + a} \right)}}}.
\end{align}
Combining \eqref{c5} and \eqref{d12}, we obtain the result.\hfill$\square$
\begin{thm}\label{thmd-3} For $a\in \mathbb{C}\setminus\Z^-$ and $k\in \N\setminus\{1\}$. Then
\begin{align}\label{d13}W_k^{\left( {a,a} \right)}\left( {m,1} \right) = \sum\limits_{r = 1}^k {{{\left( { - 1} \right)}^{r + 1}}\left( {\begin{array}{*{20}{c}}
   k  \\
   r  \\
 \end{array} } \right)\left\{ \begin{gathered}
 \sum_{i=1}^{m-1} \frac{(-1)^{i-1}}{a^i}\zeta(m+1-i)-(-1)^m \frac{H_a}{a^m}  \hfill \\ + \sum\limits_{j = 1}^{m - 1} {{{\left( { - 1} \right)}^{j - 1}}\zeta \left( {m + 1 - j} \right)H_{r - 1}^{\left( j \right)}\left( a \right)}  \hfill \\
   + {\left( { - 1} \right)^{m - 1}}{H_a}H_{r - 1}^{\left( m \right)}\left( a \right) + {\left( { - 1} \right)^{m - 1}}\sum\limits_{i = 1}^{r - 1} {\frac{{H_i^{\left( 1 \right)}\left( a \right)}}
{{{{\left( {i + a} \right)}^m}}}}  \hfill \\
\end{gathered}  \right\}} .
\end{align}
\end{thm}
\pf By a similar argument as in the proof of Theorem \ref{thmd-1}, by using \eqref{d8} and letting $f(n)=H^{(m)}_n$ with $p=1$, we have
\begin{align}\label{d14}
W_k^{\left( {a,a} \right)}\left( {m,1} \right) =& \sum\limits_{n = 1}^\infty  {\frac{{H_n^{\left( m \right)}}}
{{\left( {n + a} \right)\left( {\begin{array}{*{20}{c}}
   {n + k + a}  \\
   k  \\
 \end{array} } \right)}}}\nonumber \\
 =& k\sum\limits_{r = 1}^{k - 1} {{{\left( { - 1} \right)}^{r + 1}}r\left( {\begin{array}{*{20}{c}}
   {k - 1}  \\
   r  \\
 \end{array} } \right)} \sum\limits_{n = 1}^\infty  {\frac{{H_n^{\left( m \right)}}}
{{\left( {n + a} \right)\left( {n + r + a} \right)}}}.
\end{align}
Substituting \eqref{c5} into \eqref{d14} yields the desired result.\hfill$\square$

In the same manner we also obtain the following Theorems.
\begin{thm}\label{thmd-4} For $a\geq 0$ and $k\in \N\setminus\{1\}$. Then the following identity holds:
\begin{align}\label{d15}
W_k^{\left( {a,b} \right)}\left( {1,1,0} \right) = k\sum\limits_{r = 1}^{k - 1} {{{\left( { - 1} \right)}^{r + 1}}\left( {\begin{array}{*{20}{c}}
   {k - 1}  \\
   r  \\
 \end{array} } \right)\left\{ \begin{gathered}
  \zeta \left( 2 \right)H_r^{\left( 1 \right)}\left( a \right) - {H_a}H_r^{\left( 2 \right)}\left( a \right) \hfill \\
   - \sum\limits_{i = 1}^{r - 1} {\frac{{H_i^{\left( 1 \right)}\left( {a + 1} \right)}}
{{{{\left( {i + a + 1} \right)}^2}}}}  + \sum\limits_{j = 1}^r {\frac{{H_{a + j}^2 + H_{a + j}^{\left( 2 \right)}}}
{{a + j}}}  \hfill \\
\end{gathered}  \right\}}.
\end{align}
\end{thm}
\begin{thm}\label{thmd-5} For $a>0$ and $k\in \N$. Then
\begin{align}\label{d16}
W_k^{\left( {a,a} \right)}\left( {1,1,1} \right) = \sum\limits_{r = 1}^k {{{\left( { - 1} \right)}^{r + 1}}\left( {\begin{array}{*{20}{c}}
   k  \\
   r  \\
 \end{array} } \right)\left\{ \begin{gathered}
  \zeta \left( 2 \right)H_r^{\left( 1 \right)}\left( {a - 1} \right) - {H_a}H_r^{\left( 2 \right)}\left( {a - 1} \right) \hfill \\
   - \sum\limits_{i = 1}^{r - 1} {\frac{{H_i^{\left( 1 \right)}\left( a \right)}}
{{{{\left( {i + a} \right)}^2}}}}  + \sum\limits_{j = 1}^r {\frac{{H_{a + j - 1}^2 + H_{a + j - 1}^{\left( 2 \right)}}}
{{a + j - 1}}}  \hfill \\
\end{gathered}  \right\}}.
\end{align}
\end{thm}

Therefore, from equation \eqref{b6} and Theorems \ref{thmd-1}-\ref{thmd-5}, we know that all Euler type sums of the form
\[W_k^{\left( {a,b} \right)}\left( {1,p} \right),\ W_k^{\left( {a,b} \right)}\left( {m,0} \right),\ W_k^{\left( {a,a} \right)}\left( {m,1} \right),\ W_k^{\left( {a,b} \right)}\left( {1,1,0} \right),\ W_k^{\left( {a,a} \right)}\left( {1,1,1} \right)\]
can be expresses as a rational linear combination of products of parametric harmonic numbers, shifted harmonic numbers and zeta values. By using formulas \eqref{c20}, \eqref{c21}, \eqref{d7} and \eqref{d8}, then the Euler type sums
\[W_k^{\left( {a,a} \right)}\left( {1,2,p}\right)\quad {\rm and}\quad W_k^{\left( {a,a} \right)}\left( {1,1,1,p}\right)\]
are reducible to parametric harmonic numbers, shifted harmonic numbers, zeta values and linear sums, where $p=0$ and $1$.

{\bf Acknowledgments.} \textcolor{red}{All authors thank the anonymous referee for many invaluable comments and suggestions which have improved the paper greatly.} Ce Xu is supported by the National Natural Science Foundation of China (Grant No. 12101008), the Natural Science Foundation of Anhui Province (Grant No. 2108085QA01) and the University Natural Science Research Project of Anhui Province (Grant No. KJ2020A0057).
 \end{document}